\documentclass[12pt]{article}

\usepackage{graphicx}
\usepackage{amsmath,amssymb}
\usepackage{epstopdf}
\usepackage{psfrag}
\usepackage{float}
\usepackage[small,nohug,heads=littlevee]{diagrams} \diagramstyle[labelstyle=\scriptstyle] 
\numberwithin{equation}{section}

\begin{document}
\date{}

\begin{flushright}{YITP-SB-07-40}\end{flushright}

\begin{center}
{\LARGE Elliptic constructions of hyperk\"ahler metrics II: \\[5pt] The quantum mechanics of a Swann bundle}
\vskip40pt
{\bf Radu A. Iona\c{s}}
\end{center}
\vskip10pt
\centerline{\it C.N.Yang Institute for Theoretical Physics, Stony Brook University}
\centerline{\it Stony Brook, NY 11794-3840, USA }
\centerline{\tt ionas@max2.physics.sunysb.edu}
\vskip60pt

\begin{abstract}
\noindent The generalized Legendre transform method of Lindstr\"om and Ro\v{c}ek yields hyperk\"ahler metrics from holomorphic functions. Its main ingredients are sections of ${\cal O}(2j)$ bundles over the twistor space satisfying a reality condition with respect to antipodal conjugation on the hyperk\"ahler sphere of complex structures. Formally, the structure of the real ${\cal O}(2j)$ sections is identical to that of quantum-mechanical wave functions describing the states of a particle with spin $j$ in the spin coherent representation. We analyze these sections and their $SO(3)$ invariants and illustrate our findings with two Swann bundle constructions.
\end{abstract}
\vspace*{40pt}

\newpage

\tableofcontents

\setcounter{section}{-1}

\section{Introduction}

This paper continues the program initiated in \cite{MM1}, concerning the generalized Legendre transform method of constructing hyperk\"ahler metrics \cite{Lindstrom:1983rt,Hitchin:1986ea}. At the center of this approach stands a single holomorphic function of certain real sections of ${\cal O}(2j)$ bundles over the twistor space of the hyperk\"ahler variety that encodes all the metric information. In the first part of these notes, we explore a quantum-mechanical analogy along some ideas of Penrose to gain insight into the structure of these ${\cal O}(2j)$ sections. The second part deals with applications.

A direction of applications is the construction of Swann bundle metrics \cite{MR1096180}. Swann bundles, known also as hyperk\"ahler cones,  are hyperk\"ahler varieties with an additional $\mathbb{H}^*$-action, whose real component acts homothetically while the three purely imaginary components act isometrically and rotate the hyperk\"ahler complex structures. The complex structures have furthermore the distinctive feature that they  can be derived from the same K\"ahler potential, defined up to the addition of a constant. The importance of this class of hyperk\"ahler varieties stems from the fact that any quaternionic-K\"ahler manifold possesses a canonical Swann bundle from which it can be obtained through a quotient construction. They thus provide a holomorphic environment for the description of the generally non-holomorphic  quaternionic-K\"ahler manifolds. In \cite{MA1}, drawing on an observation of \cite{deWit:2001dj}, we proved a  criterion for a generalized Legendre transform approach to yield a hyperk\"ahler variety with a Swann bundle structure: the aforementioned holomorphic function must be a section (sometimes an affine section) of an ${\cal O}(2)$ bundle over the twistor space.

A different motivation for our choice of applications derived from a question that emerged during the study of the Atiyah-Hitchin metric in the generalized Legendre transform frame. To formulate it, let us recall a few facts from \cite{MM1}. To each ${\cal O}(4)$ real section one associates canonically a quartic plane curve having, among other features, an orthogonal period lattice. On the other hand, on the parameter space of these sections there exists an $SO(3)$ action induced by the automorphisms of the Riemann sphere that preserve the real structure of the sections. Under this action, the parameters transform in the spin-$2$ representation of $SO(3)$. The two Casimir invariants of the action can be taken to be the periods of the ${\cal O}(4)$ curve. In fact, from a practical point of view it is more convenient to consider the reciprocals of the periods instead, which we denote, up to a numerical factor, by $r_2$ and $r'_2$. The orthogonality of the lattice guarantees that we can choose the numerical factor such that $r_2,r'_2 > 0$. So these look like radial coordinates. In terms of them, the elliptic {\it nome} and complementary {\it nome} take the forms $q = \exp(-\pi^2r_2/r'_2)$ and $q' = \exp(-r'_2/r_2)$ respectively, hence the two asymptotic regions $r_2>>r'_2$ and $r_2<<r'_2$ can be analyzed perturbatively by performing series expansions in $q$ respectively $q'$. In the Atiyah-Hitchin case, the defining equation of the manifold in the generalized Legendre approach turns out to be simply $r_2 = h = \mbox{const}$. The remaining radius, $r'_2$, plays the role of monopole separation distance. The limits $r'_2>>h$ and $r'_2<<h$ have been studied perturbatively in \cite{Atiyah:1988jp,Gibbons:1995yw,Hanany:2000fw}.  One might ask the question, are there examples in which none of the radii is frozen? In this paper we answer this question positively and provide such an example.

\section{Spherical representations and invariants of ${\cal O}(2j)$ multiplets} \label{SEC:rot-inv}

${\cal O}(2j)$ multiplets are by definition ${\cal O}(2j)$ sections $\eta^{(2j)} = \eta_{A_1\cdots A_{2j}}\pi^{A_1} \cdots \pi^{A_{2j}}$ over the twistor space $Z$ (we can define sections of ${\cal O}(k)$ over $Z$ by pulling back from $\mathbb{CP}^1$) which satisfy a reality condition with respect to antipodal conjugation. Letting $\zeta = \pi^2/\pi^1$ be the standard inhomogeneous coordinate on $\mathbb{CP}^1$, we define the local 'tropical' form $\eta^{(2j)}(\zeta)$ by
\begin{equation}
\eta_{A_1\cdots A_{2j}}\pi^{A_1} \cdots \pi^{A_{2j}} = (\pi^1\pi^2)^j \eta^{(2j)}(\zeta) 
\end{equation}
The reality constraint then reads
\begin{equation}
\eta^{(2j)}(-\frac{1}{\bar{\zeta}}) = \overline{\eta^{(2j)}(\zeta)}\label{reality1}
\end{equation}
A generic ${\cal O}(2j)$ multiplet can be written in either one of the following two equivalent local forms
\begin{eqnarray}
\eta^{(2j)}(\zeta) &\hspace{-7pt} = \hspace{-7pt}& \sum_{m=-j}^{j} \!\left(\!\! \begin{array}{c} 2j \\ j\!+\!m \end{array} \!\!\right)^{\!\!1/2}\! \bar{\psi}^{\,j}_m \, \zeta^m \label{Maj_coeffs} \\
&\hspace{-7pt} = \hspace{-7pt}& \frac{\varrho}{\zeta^j} \prod_{l=1}^j \frac{(\zeta-\alpha_l)(\bar{\alpha}_l\zeta+1)}{1+|\alpha_l|^2 } \label{Maj_roots}
\end{eqnarray}
The requirement (\ref{reality1}) translates into the condition $\psi^{\,j}_{-m} = (-)^m \bar{\psi}^{\,j}_m$ on the coefficients in the first line and the condition $\varrho \in \mathbb{R}$ as well as into the antipodal pairing of the roots in the second line. 

The reality constraint is preserved only by the $\mathbb{R}^* \times SU(2)$ subgroup of the $PSL(2,\mathbb{C})$ group of automorphisms of the Riemann sphere which act through birational transformations on the inhomogeneous coordinate $\zeta$. Under the M\"{o}bius action of an element $R$ of its $SU(2)$ component 
\begin{eqnarray}
\zeta & \stackrel{R}{\longrightarrow} & \frac{\hspace{5pt} a\, \zeta+b}{-\bar{b}\, \zeta+\bar{a}} \label{Moebius}
\end{eqnarray}
with $a,b \in \mathbb{C}$ such that $|a|^2+|b|^2=1$, the roots of $\eta^{(2j)}(\zeta)$ transform in the same way, {\it i.e.},
\begin{eqnarray}
\alpha_l & \stackrel{R}{\longrightarrow} & \frac{\hspace{5pt} a\, \alpha_l+b}{-\bar{b}\, \alpha_l+\bar{a}} 
\end{eqnarray}
whereas the scale factor $\varrho$ remains inert. Thus, the root system consists of $j$ antipodal pairs of points on the Riemann sphere that rotate together rigidly - a {\it constellation} in the language of \cite{Bacry:1974hc,MR0479078}. On the other hand, under (\ref{Moebius}) the coefficients $\psi^{\,j}_m$ transform according to Wigner's D-function realization of the spin-$j$ unitary irrep of $SO(3)$ - the double-cover of $SU(2)$,
\begin{eqnarray}
\hspace{63pt} \psi^{\,j}_m & \stackrel{R}{\longrightarrow} & \sum_{m'=-j}^{j} D^{\,j}_{mm'}(\phi,\theta,\psi)\psi^{\,j}_{m'} \label{Wigner}
\end{eqnarray}
where the Euler angles are related to the Cayley-Klein parameters of $R$ by
\begin{equation}
a = \cos\! \frac{\theta}{2}\, e^{\frac{i}{2}(\phi+\psi)} \qquad\qquad b = \sin\! \frac{\theta}{2}\, e^{\frac{i}{2}(\phi-\psi)}
\end{equation}

\subsection{Quantum spin coherent states} \label{coherent_states}

Polynomials of the type $\zeta^j \eta^{(2j)}(\zeta)$ occur in the context of Quantum Mechanics in the guise of (unnormalized) spin-$j$ wave functions. In the form (\ref{Maj_coeffs}) they are known as being in the spin coherent state representation \cite{Radcliffe:1971}, whereas in the form (\ref{Maj_roots}) as being in Majorana's  {\it stellar representation} \cite{Majorana:1932ga}. For this latter reason we shall refer to them in these notes as Majorana polynomials.

The quantum states of a particle with spin $j$ are commonly described as linear superpositions of $2j+1$ spherical harmonics, {\it i.e.},
\begin{equation}
| \psi \rangle = \sum_{m=-j}^{j} \psi^j_m \, | j m \rangle \label{ket}
\end{equation}
The spherical harmonics are simultaneous eigenvalues of the Casimir operator $J^2$ and of the operator $J_z$ corresponding to  the projection of the angular momentum along a preferential axis and form an orthonormal basis in the Hilbert space of states that transforms according to the spin-$j$ unitary irrep of $SU(2)$. Under such a transformation, the linear coefficients of (\ref{ket}) transform just as in (\ref{Wigner}).
Corresponding to any $\zeta \in \mathbb{C} \cup \{\infty\}$ one defines in this basis a {\it spin coherent state} by \cite{Radcliffe:1971}
\begin{equation}
| \zeta \rangle = \frac{1}{(1+|\zeta|^2)^j} \sum_{m=-j}^{j} \!\left(\!\! \begin{array}{c} 2j \\ j\!+\!m \end{array} \!\!\right)^{\!\!1/2}\!  \zeta^{j+m} | j m \rangle 
\end{equation}
The set of spin coherent states forms an overcomplete basis in the Hilbert space of states. The wave function $\psi(\zeta) = \langle\psi|\zeta\rangle$ is said to be the spin coherent state representation of $|\psi\rangle$, and is equal, up to a non-holomorphic normalization factor, to a Majorana polynomial holomorphic in $\zeta$, {\it i.e.}
\begin{equation}
\langle \psi | \zeta \rangle = \frac{1}{(1+|\zeta|^2)^j} \zeta^j \eta^{(2j)}(\zeta) \label{w_function}
\end{equation}
with $\eta^{(2j)}(\zeta)$ expressed as in (\ref{Maj_coeffs}). In particular, the spin coherent state representation of a purely spin coherent state labeled by the complex number $\alpha$ takes the form
\begin{equation}
\langle \alpha | \zeta \rangle = \left[\frac{1+\bar{\alpha}\zeta}{\sqrt{(1+|\zeta|^2)(1+|\alpha|^2)}} \right]^{2j} \label{coh_st}
\end{equation}
Corresponding to the factorization (\ref{Maj_roots}), the spin-$j$ wave function (\ref{w_function}) decomposes, up to a quantum-mechanically irrelevant phase factor, into a product of $2j$ spin-$1/2$ coherent wave functions
\begin{equation}
\langle \psi | \zeta \rangle \sim \varrho \prod_{l=1}^j \langle-\frac{1}{\bar{\alpha}_l}|\zeta\rangle_{1/2} \langle\alpha_l|\zeta\rangle_{1/2}
\end{equation}
Similarly, the coherent wave function (\ref{coh_st}) can be written as
\begin{equation}
\langle\alpha|\zeta\rangle = [ \langle\alpha|\zeta\rangle_{1/2}]^{2j}
\end{equation}
A very intuitively appealing picture emerges: a quantum state with spin $j$ appears to be described by a set of $2j$ elementary "spins 1/2" with the origins at the center of a Bloch sphere, pointing out in the directions marked by a constellation of $2j$ dots on surface of the sphere corresponding to the roots of the wave function polynomial. In particular, a spin state is real in the sense of (\ref{reality1}) when all elementary spins come in oppositely oriented pairs and is coherent when all elementary spins point in the same direction.
Clearly, in the spin coherent state representation the rotational structure is preserved manifestly and no preferential axis needs to be chosen.

These elementary spins correspond essentially to Penrose's notion of {\it principal spinors}, as defined {\it e.g.} in  \cite{Penrose:1986ca}, see also \cite{MR1865778}. Penrose frames the above result in the following language: any nonvanishing totally symmetric $\eta_{A_1\cdots A_{2j}}$ admits a canonical decomposition 
\begin{equation}
\eta_{A_1\cdots A_{2j}} = \chi^{(1)}_{(A_1}\chi^{(2)}_{A_2\phantom{)}} \!\cdots \chi^{(2j)}_{A_{2j})} \label{Penrose_fact}
\end{equation}
in terms of a set of $2j$ commutative spinors $\chi^{(k)}_A$, uniquely defined up to proportionality and reordering. 

The properties of quantum spin-1/2 coherent states are especially fit for use in spherical geometry, and we will exploit this feature later on. The overlap between two spin-1/2 coherent states corresponding to $\alpha, \beta \in \mathbb{C} \cup \{\infty\} \simeq S^2$ is\footnote{To avoid cluttering, for the remainder of this section we drop the index 1/2 from the notation of spin-1/2 coherent wave functions.} 
\begin{equation}
\langle \alpha | \beta \rangle = \frac{1+\bar{\alpha}\beta}{\sqrt{(1+|\alpha|^2)(1+|\beta|^2)}} 
\end{equation}
Note that this formula implies that the overlap between states corresponding to pairs of antipodally-conjugated points is zero. The norms
\begin{equation}
|\langle \alpha | \beta \rangle| = k_{\alpha\beta} \hspace{40pt} |\langle -\frac{1}{\bar{\alpha}} | \beta \rangle| = k'_{\alpha\beta} \label{coh_norms}
\end{equation}
are related to the geodesic distance on the sphere between $\alpha$ and $\beta$, see equations~(\ref{FS}) and (\ref{chordal_dist}) below. On the other hand, the phases of cyclic sequences of spin-1/2 coherent states have an area interpretation, namely,
\begin{equation}
\langle \alpha_1 | \alpha_2 \rangle \langle \alpha_2 | \alpha_3 \rangle \cdots \langle \alpha_{n-1} | \alpha_n \rangle \langle \alpha_n | \alpha_1 \rangle = k_{\alpha_1\alpha_2} k_{\alpha_2\alpha_3} \cdots k_{\alpha_{n-1}\alpha_n} k_{\alpha_n\alpha_1} e^{i A_{\rm polygon}/2} \label{cyclic_coh}
\end{equation}
where $A_{\rm polygon}$ is the area of the spherical polygon with vertices at the points $\alpha_1$ $\cdots$  $\alpha_n$. The factor 1/2 in front of the area makes the ambiguity in the choice of what one calls the "inside" and the "outside" of the polygon irrelevant.

\subsection{Rotational invariants}

For reasons to become clear later on, we are interested in constructing $SU(2)$ invariant quantities associated to one multiplet or a system of multiplets. To this purpose, we develop several  complementary approaches. 

Our first approach is a natural by-product of the geometric picture detailed above. Consider a multiplet or a set of multiplets for which we want to compute invariants and the corresponding constellation of roots on the Riemann sphere endowed with the $SU(2)$-invariant metric of Fubini and Study. Given two such roots $\alpha$ and $\beta$ from the same or from two different multiplets, the Fubini-Study distance between them is given by
\begin{equation}
\delta_{\alpha\beta} = 2 \arccos k_{\alpha\beta} = 2 \arcsin k'_{\alpha\beta} \label{FS}
\end{equation}
with the chordal distance and radius expressed in terms of the roots as follows
\begin{equation}
k_{\alpha\beta} = \frac{|1+\bar{\alpha}\beta|}{\sqrt{(1+|\alpha|^2)(1+|\beta|^2)}} \hspace{30pt} {\rm and} \hspace{30pt}  k'_{\alpha\beta} = \frac{|\alpha-\beta|}{\sqrt{(1+|\alpha|^2)(1+|\beta|^2)}} \label{chordal_dist}
\end{equation}
One has $k^2_{\alpha\beta} + k_{\alpha\beta}^{\prime2}  = 1$ and thus $0 < k_{\alpha\beta}, k_{\alpha\beta}' < 1$. We can then use invariant Fubini-Study distances as building blocks to construct proper invariants by considering combinations of them subject to the additional condition that they are symmetric at the permutation of the roots of each of the multiplets involved. 

A second approach involves constructing invariant Penrose transforms. It is based on the following result: let
\begin{equation}
{\cal I} = \oint_{\Gamma} \frac{d\zeta}{\zeta} G(\eta^{(2j)}(\zeta)) \label{Penrs-int}
\end{equation}
be a contour integral, with $G$ a meromorphic function possibly with branch cuts, depending on one or several multiplets denoted here collectively by $\eta^{(2j)}$ and $\Gamma$ an integration contour that yields either a real or a purely imaginary ${\cal I}$. If   \\
\phantom{lo} 1)\hspace{1pt} $G$ does not depend explicitly on $\zeta$ other than through $\eta^{(2j)}(\zeta)$, and \\
\phantom{lo} 2) $G$ scales, modulo terms that vanish under the contour integral, with weight $-1$ when each $\eta^{(2j)}(\zeta)$ is scaled with weight $j$, \\
\noindent then ${\cal I}$ can be cast in the following manifestly $SU(2)$-invariant form
\begin{equation}
{\cal I} = \oint_{\Gamma} \pi_Ad\pi^A G(\eta_{A_1\cdots A_{2j}}\pi^{A_1} \!\cdots\, \pi^{A_{2j}})
\end{equation}

Alternatively, one can use the spherical tensor properties of the coefficients $\psi^{\,j}_{m}$ of the Majorana polynomials to construct spherical scalars by invariantly coupling two such tensors, three, a.s.o., {\it i.e.},
\begin{equation}
\sum_m \bar{\psi}^{\,j}_m \psi^{\,j}_m \ ,
\quad 
\sum_{m_1,m_2,m_3} \!\left(\!\!\begin{array}{ccc} j_1&j_2&j_3 \\ m_1&m_2&m_3 \end{array} \!\!\right) \psi^{\,j_1}_{m_1} \psi^{\,j_2}_{m_2} \psi^{\,j_3}_{m_3} \ ,
\quad \cdots \label{spherical_invs}
\end{equation}
The coupling factors in the second expression are Wigner $3j$-symbols. The formulas become increasingly more complex with the number of angular momenta coupled.

A more uniform approach that leads to equivalent results is to form spherical scalars by completely contracting indices of various combinations of symmetric tensors $\eta_{A_1\cdots A_{2j}}$ corresponding to the set of multiplets one is interested in computing invariants for. At first sight it may look like there exists an infinite number of such scalar configurations, but Penrose's canonical decomposition (\ref{Penrose_fact})  implies that only a finite number of them are in fact independent. There is a nice way to depict these scalar combinations graphically by representing {\it e.g.} $\eta_{A_1\cdots A_{2j}}$ as a vertex with $2j$ outgoing lines and $\eta^{A_1\cdots A_{2j}}$ as a similar vertex but with $2j$ incoming lines. The resulting graphs can then be easily manipulated and related to each other by using diagrammatic identities such as \\ [-5pt]
\begin{figure}[H]
\centering
\scalebox{0.55}{\includegraphics{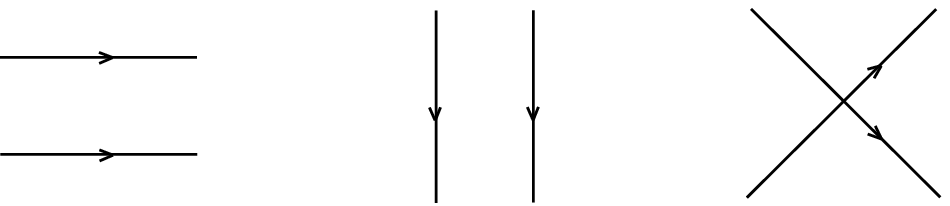}}
\put(-105,13){=}
\put(-45,13){-}
\end{figure} 
\noindent which expresses the $\epsilon$-symbol identity $\epsilon^{AB}\epsilon_{CD}=\delta^A{}_C \delta^B{}_D - \delta^A{}_D \delta^B{}_C$. Graphs with legs starting and ending on the same vertex vanish, reflecting the fact that $\eta_{A_1\cdots A_{2j}}$ is totally symmetric in its indices and hence yields zero when two of these are contracted with an $\epsilon$-symbol. Reversing the orientation of a leg changes the sign of the graph.

\subsection{Invariants of ${\cal O}(2)$ multiplets} \label{O2-invs}

A generic ${\cal O}(2)$ multiplet can be written locally in either one of the following two forms
\begin{eqnarray}
\eta^{(2)}(\zeta)	&\hspace{-7pt} = \hspace{-7pt}& \frac{\bar{z}_1}{\zeta}+x_1-z_1\zeta \nonumber \\
			&\hspace{-7pt} = \hspace{-7pt}& \frac{\sigma}{\zeta}\frac{(\zeta-\gamma)(\bar{\gamma}\zeta+1)}{1+|\gamma|^2}
\end{eqnarray}
The coefficients can be expressesed in terms  of the roots and scale factor $\sigma$ explicitly. Conversely, in order to express the roots in terms of the coefficients one has to solve a quadratic equation.

\begin{figure}[ht]
\centering
\scalebox{0.55}{\includegraphics{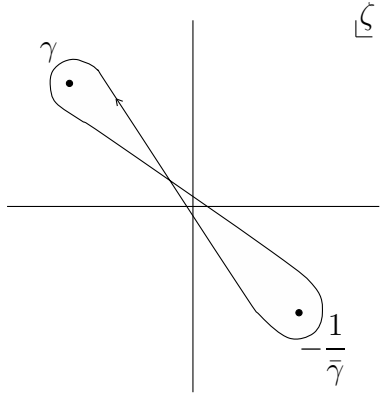}}
\put(-129,128){$\gamma$}
\put(-31,14){$\displaystyle{-\frac{1}{\bar{\gamma}}}$}
\put(-8,140){$\zeta$}
\caption{The contour $\Gamma$}
\label{int_inv_1}
\end{figure} 
To a real ${\cal O}(2)$ section one can associate only one independent invariant, namely $\sigma$. An invariant integral is
\begin{equation}
\oint_{\Gamma} \frac{d\zeta}{\zeta} \frac{1}{\eta^{(2)}(\zeta)} = \frac{2}{\sigma} \label{O2-radius-oint}
\end{equation}
with the contour $\Gamma$ depicted in \figurename~\ref{int_inv_1}.

Alternatively, as suggested above, one can construct the scalar 
\begin{figure}[H]
\centering
\hspace{120pt}
\scalebox{0.4}{\includegraphics{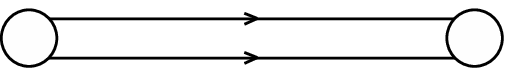}}
\put(-188,1){$g_{\sigma^2} = -2\, \eta_{AB}\eta^{AB} = -2 \times$}
\end{figure}
\vspace{-4pt}
\noindent The factor $-2$ has been inserted for convenience. A short calculation yields that
\begin{equation}
g_{\sigma^2} = 4 |z_1|^2+x_1^2  = \sigma^2
\end{equation}
As one can easily check, all other invariants can be deconstructed down to $g_{\sigma^2}$, {\it e.g.},
\begin{figure}[H]
\centering
\hspace{-40pt}
\scalebox{0.4}{\includegraphics{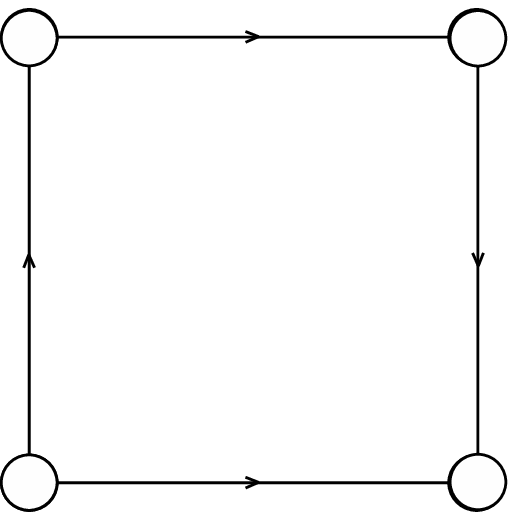}}
\put(2,27){$= \displaystyle{\frac{1}{8} g_{\sigma^2}^{\,2}}$}
\end{figure} 
\begin{figure}[H]
\centering
\hspace{-40pt}
\scalebox{0.4}{\includegraphics{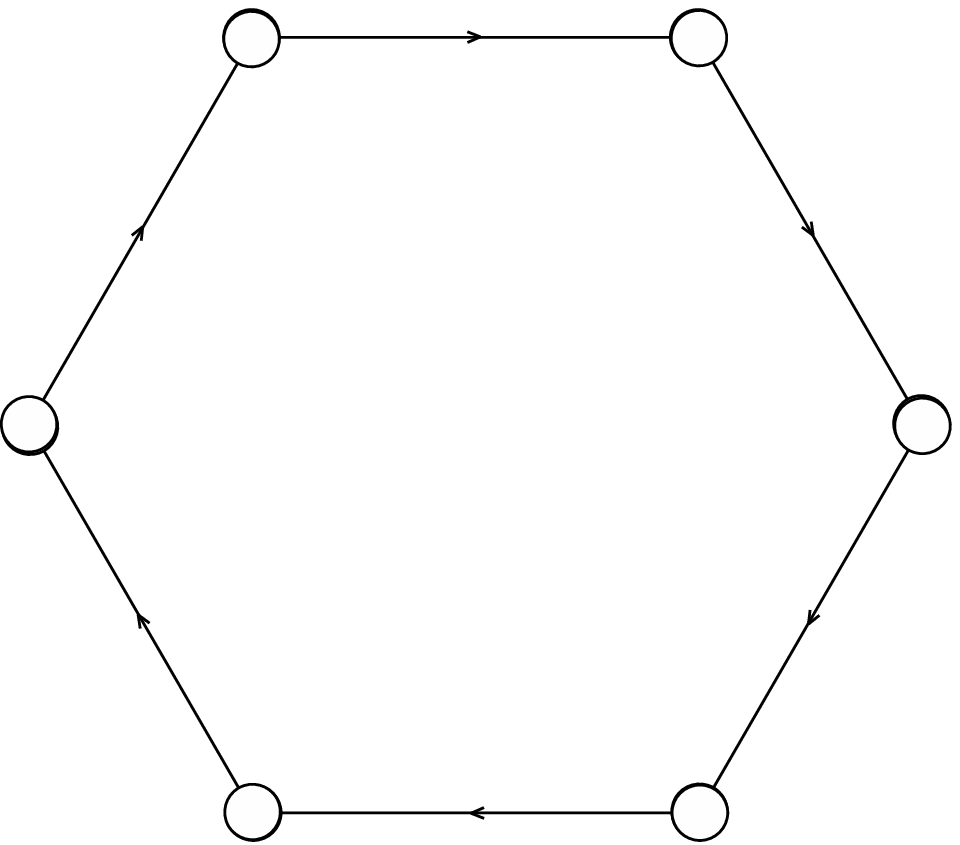}}
\put(3,45){$= - \displaystyle{\frac{1}{32} g_{\sigma^2}^{\,3}}$}
\end{figure} 
\noindent a.s.o. Polygons with $2k+1$ sides vanish identically. Polygons with $2k$ sides yield the $k$-th power of $g_{\sigma^2}$, times a numerical factor.

The coefficients $z_1$, $x_1$, $\bar{z}_1$ form an $SO(3)$ vector multiplet. This can be cast in an Euclidean basis by the linear transformation
\begin{equation}
z_1 = \frac{1}{2}(x+iy) \qquad\qquad 
x_1 = z \qquad\qquad
\bar{z}_1 = \frac{1}{2}(x-iy)
\end{equation}
We shall use the notation $\vec{r}_1$ for the $\mathbb{R}^3$ vector with components $x$, $y$, $z$. Clearly,
\begin{equation}
|\vec{r}_1| = \sqrt{x^2+y^2+z^2} = \sigma
\end{equation}
{\it i.e.}, $\sigma$ represents the invariant length of the vector associated in this manner with the ${\cal O}(2)$ multiplet.

\subsection{Invariants of ${\cal O}(4)$ multiplets} \label{O4-invs}

A generic ${\cal O}(4)$ multiplet can be written locally in either one of the following two forms
\begin{eqnarray}
\eta^{(4)}(\zeta) &\hspace{-7pt} = \hspace{-7pt}& \frac{\bar{z}_2}{\zeta^2}+\frac{\bar{v}_2}{\zeta}+x_2-v_2\zeta+z_2\zeta^2 \nonumber \\
&\hspace{-7pt} = \hspace{-7pt}&  \frac{\rho}{\zeta^2} \frac{(\zeta-\alpha)(\bar{\alpha}\zeta+1)}{1+|\alpha|^2}\frac{(\zeta-\beta)(\bar{\beta}\zeta+1)}{1+|\beta|^2} \label{eta_four}
\end{eqnarray}
The coefficients can be expressed directly in terms of the roots, but conversely, expressing the roots explicitly in terms of the coefficients involves solving a quartic equation, an impractical approach.

\begin{figure}[tbh]
\centering
\scalebox{0.75}{\includegraphics{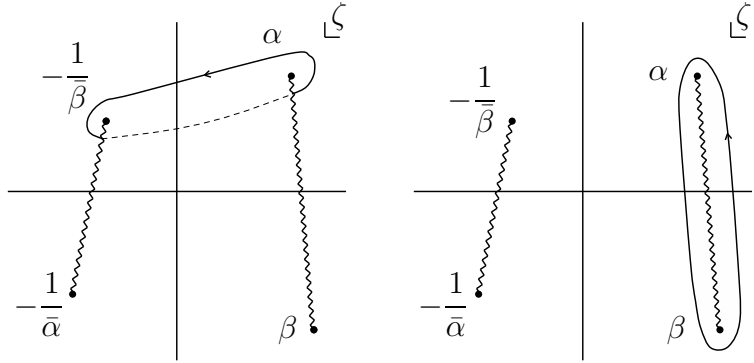}}
\put(-40,107){$\alpha$}
\put(-127,18){$\displaystyle{-\frac{1}{\bar{\alpha}}}$}
\put(-33,9){$\beta$}
\put(-116,96){$\displaystyle{-\frac{1}{\bar{\beta}}}$}
\put(-7,128){$\zeta$}
\put(-186,121){$\alpha$}
\put(-280,18){$\displaystyle{-\frac{1}{\bar{\alpha}}}$}
\put(-180,9){$\beta$}
\put(-270,105){$\displaystyle{-\frac{1}{\bar{\beta}}}$}
\put(-160,128){$\zeta$}
\caption{Integration contours $\Gamma_a$ (left) and $\Gamma_b$ (right).}
\label{int_inv_2}
\end{figure} 
An invariant integral with the required homogeneity property is
\begin{equation}
{\cal I}(\Gamma) = \oint_{\Gamma} \frac{d\zeta}{\zeta} \frac{1}{\sqrt{\eta^{(4)}(\zeta)}}
\end{equation}
The two generators of the canonical homology basis for the closed contours $\Gamma$ are depicted in \figurename~\ref{int_inv_2}. They correspond to the $a$ and $b$-cycles of the ${\cal O}(4)$ curve associated to the multiplet \cite{MM1}. The integrals over these two contours are precisely the period integrals of the ${\cal O}(4)$ curve. In \cite{MM1} we have shown that they can be expressed in terms of the complete elliptic integrals of modulus $k_{\alpha\beta}$ respectively complementary modulus $k'_{\alpha\beta}$ as follows
\begin{equation}
{\cal I}(\Gamma_a) = \frac{2}{\sqrt{\rho}} K(k_{\alpha\beta}) \stackrel{\rm def}{=} \frac{2}{r_2} \qquad\qquad
{\cal I}(\Gamma_b) = \frac{2}{\sqrt{\rho}} iK(k'_{\alpha\beta}) \stackrel{\rm def}{=} \frac{2\pi i}{r'_2} \label{O4-radii-oint}
\end{equation}
The second set of equalities are definitions inspired by and analogous to (\ref{O2-radius-oint}). The $\pi$-factor in the r.h.s. has been chosen for later convenience but is otherwise irrelevant. Since $0<k_{\alpha\beta},k'_{\alpha\beta}<1$, the elliptic integrals are real and so 
\begin{equation}
r_2, r'_2 > 0
\end{equation}
Based on this and their rotational invariance property, we shall refer to $r_2$ and $r'_2$ as "${\cal O}(4)$ radii". If instead we describe the ${\cal O}(4)$ curve in the Weierstrass framework, the complete elliptic integrals are customarily replaced by the Weierstrass half-periods $\omega$ and $\omega'$, in terms of which the definitions (\ref{O4-radii-oint}) read \cite{MM1}
\begin{equation}
r_2=\frac{1}{2\omega} \qquad \mbox{and} \qquad r'_2=\frac{i\pi}{2\omega'} \label{r&r'}
\end{equation}

We can also construct ${\cal O}(4)$ rotational invariants by completely contracting the indices of products of $\eta_{ABCD}$ tensors. With the diagrammatic conventions introduced above, let for instance
\begin{figure}[H]
\centering
\hspace{40pt}
\scalebox{0.4}{\includegraphics{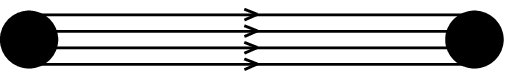}}
\put(-105,1){$\displaystyle{g_{\rho^2} = 2 \times}$}
\end{figure} 
\begin{figure}[H]
\centering
\hspace{40pt}
\scalebox{0.4}{\includegraphics{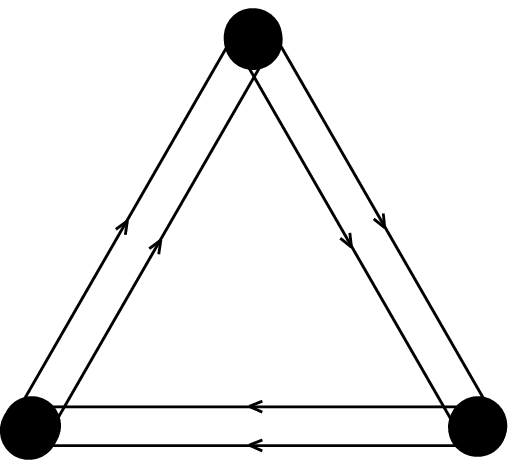}}
\put(-105,27){$\displaystyle{g_{\rho^3} = \frac{8}{3} \times}$}
\end{figure} 
The multiplicative factors have been inserted for convenience. A straightforward calculation yields the following Majorana-coefficient expressions
\begin{eqnarray}
g_{\rho^2} &\hspace{-7pt} = \hspace{-7pt}& 4 |z_2|^2 + |v_2|^2 + \frac{1}{3} x_2^2 \\
g_{\rho^3} &\hspace{-7pt} = \hspace{-7pt}&  \frac{8}{3} |z_2|^2 x_2 - \frac{1}{3} |v_2|^2 x_2 - \frac{2}{27} x_2^3 - z_2 \bar{v}_2^2 - \bar{z}_2 v_2^2 
\end{eqnarray}
$g_{\rho^2}$ and $g_{\rho^3}$ are essentially the only independent invariants that one can construct in this manner. All higher order spherical scalars break down ultimately into these two basic components. For example,
\begin{figure}[H]
\centering
\hspace{-60pt}
\scalebox{0.4}{\includegraphics{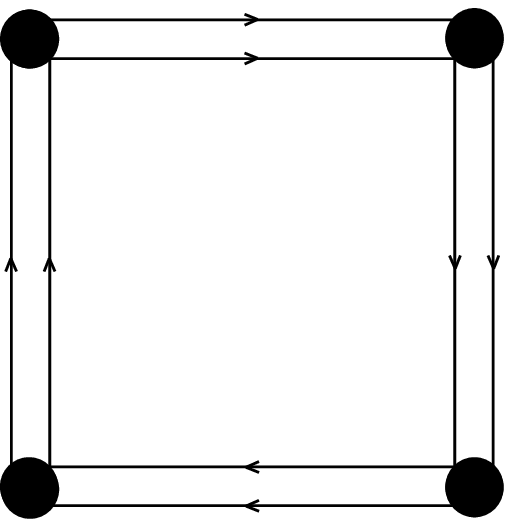}}
\put(4,27){$=\displaystyle{\frac{1}{8} g_{\rho^2}^{\,2}}$}
\end{figure} 
\begin{figure}[H]
\centering
\hspace{-60pt}
\scalebox{0.4}{\includegraphics{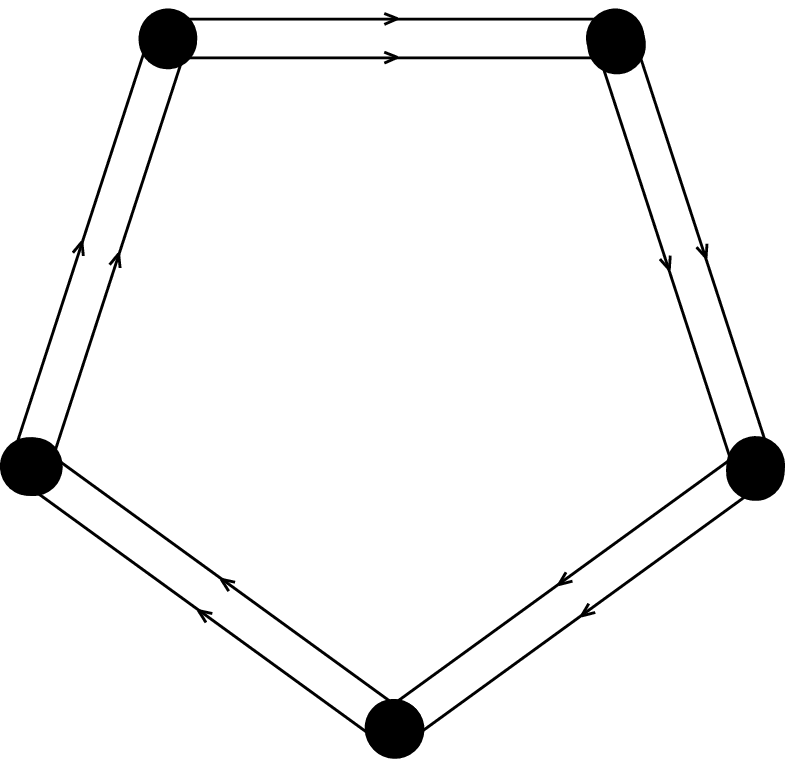}}
\put(4,42){$=\displaystyle{\frac{5}{32} g_{\rho^2} g_{\rho^3}}$}
\end{figure} 
\begin{figure}[H]
\centering
\hspace{-60pt}
\scalebox{0.4}{\includegraphics{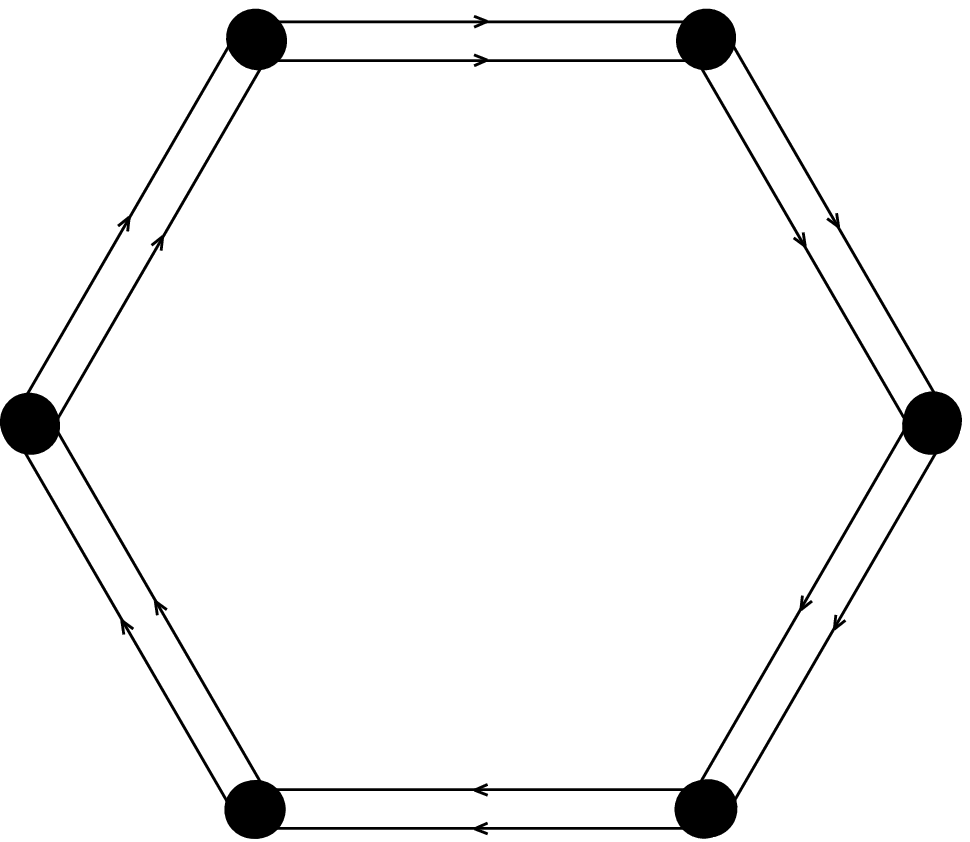}}
\put(4,46){$=\displaystyle{\frac{1}{64}(2 g_{\rho^2}^{\,3}+3g_{\rho^3}^{\,2})}$}
\end{figure} 
\noindent a.s.o. 

And yet a third pair of invariants is provided by the scale $\rho$ and the chordal distance $k_{\alpha\beta}$. Clearly though, these sets of pairs of invariants are not all independent. As a matter of fact, one can express $g_{\rho^2}$ and $g_{\rho^3}$ in terms of both $\rho$, $k_{\alpha\beta}$ and $r_2$, $r'_2$. Thus, on one hand, by resorting to the Vi\`ete relations between the coefficients and the roots of (\ref{eta_four}), one can show that
\begin{eqnarray}
g_{\rho^2} &\hspace{-7pt} = \hspace{-7pt}& -(e_1e_2+e_2e_3+e_3e_1) \\ [4pt]
g_{\rho^3} &\hspace{-7pt} = \hspace{-7pt}& e_1e_2e_3
\end{eqnarray}
with
\begin{equation}
e_1 = -\frac{\rho}{3} (k_{\alpha\beta}^2-2) \qquad\quad e_2 = \frac{\rho}{3} (2k_{\alpha\beta}^2-1) \qquad\quad e_3 = -\frac{\rho}{3}(k_{\alpha\beta}^2+1) \label{W_roots}
\end{equation}
On the other hand, in \cite{MM1} we have shown that when the ${\cal O}(4)$ curve is cast in Weierstrass form, $g_{\rho^2}$ and $g_{\rho^3}$ end up playing the role of Weierstrass modular coefficients. As is well-known in the theory of elliptic functions, the Weierstrass modular coefficients admit the following double power series representation in terms of the elliptic periods 
\begin{eqnarray}
g_{\rho^2} &\hspace{-7pt} = \hspace{-7pt}& 15\sum_{m,m'=-\infty}^{\infty} \hspace{-14pt}{}' \hspace{10pt} \frac{1}{(2m\omega+2m'\omega')^4}  \\ [2pt]
g_{\rho^3} &\hspace{-7pt} = \hspace{-7pt}& 35\sum_{m,m'=-\infty}^{\infty} \hspace{-14pt}{}' \hspace{10pt} \frac{1}{(2m\omega+2m'\omega')^6}
\end{eqnarray}
where the prime sum symbol signifies that the term with $(m,m')=(0,0)$ must be omitted. Alternatively,  each of these double series can be recast as a Lambert-type $q$-series
\begin{eqnarray}
g_{\rho^2} &\hspace{-7pt} = \hspace{-7pt}& \hspace{3pt} \frac{1}{3}\, \left(\frac{\pi}{2\omega}\right)^4 \!\left( 1+240\sum_{n=1}^{\infty} n^3\frac{q^{2n}}{1-q^{2n}} \right) \label{g2-Lambert} \\ [2pt]
g_{\rho^3} &\hspace{-7pt} = \hspace{-7pt}& \frac{2}{27} \left(\frac{\pi}{2\omega}\right)^6 \!\left( 1-504\sum_{n=1}^{\infty} n^5\frac{q^{2n}}{1-q^{2n}} \right) \label{g3-Lambert}
\end{eqnarray}
where $q=\exp(i\pi\tau)$ is the elliptic {\it nome} and $\tau=\omega'/\omega$ is the elliptic modulus. Since the Weierstrass coefficients are invariant under the modular transformation $\tau' = -1/\tau$, their $q'$-series expansions are formally identical, but with $q$ replaced by $q' =\exp(i\pi\tau')$ and $\omega$ by $\omega'$. In terms of the radii, 
\begin{equation}
 q = e^{-\pi^2r_2/r'_2} \qquad\mbox{and}\qquad q' = e^{-r'_2/r_2} 
\end{equation} 
The fact that $r_2$ and $r'_2$ are positive implies that $0<q,q'<1$, which in turn guarantees convergence. The two asymptotic regions $r_2>>r'_2$ and $r_2<<r'_2$ can be analyzed perturbatively by performing expansions in $q$ respectively $q'$.

\subsection{Mixed invariants of ${\cal O}(2)$ and ${\cal O}(4)$ multiplets} \label{mixed-invs}

Consider now the combination of an ${\cal O}(2)$ with an ${\cal O}(4)$ multiplet. An invariant integral containing both is 
\begin{equation}
{\cal I}(\Gamma) = \oint_{\Gamma} \frac{d\zeta}{\zeta} \frac{\eta^{(2)}(\zeta)}{\eta^{(4)}(\zeta)}
\end{equation}
\begin{figure}[tbh]
\centering
\scalebox{0.75}{\includegraphics{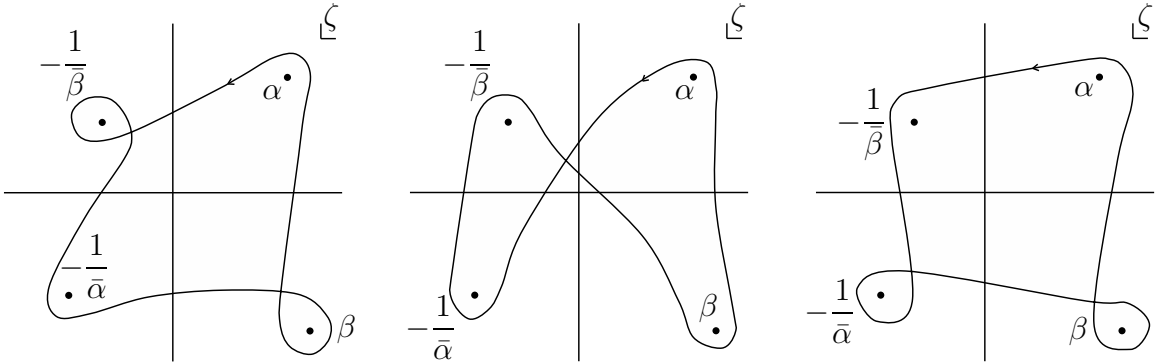}}
\put(-31,101){$\alpha$}
\put(-133,16){$\displaystyle{-\frac{1}{\bar{\alpha}}}$}
\put(-33,9){$\beta$}
\put(-121,88){$\displaystyle{-\frac{1}{\bar{\beta}}}$}
\put(-7,128){$\zeta$}
\put(-182,100){$\alpha$}
\put(-284,10){$\displaystyle{-\frac{1}{\bar{\alpha}}}$}
\put(-173,18){$\beta$}
\put(-270,110){$\displaystyle{-\frac{1}{\bar{\beta}}}$}
\put(-161,128){$\zeta$}
\put(-338,100){$\alpha$}
\put(-415,33){$\displaystyle{-\frac{1}{\bar{\alpha}}}$}
\put(-310,10){$\beta$}
\put(-423,110){$\displaystyle{-\frac{1}{\bar{\beta}}}$}
\put(-315,128){$\zeta$}
\caption{Integration contours $\Gamma_0$ (left), $\Gamma_+$ (middle) and $\Gamma_-$ (right).}
\label{int_inv_24}
\end{figure} \\ 
For the three independent contours depicted in \figurename~\ref{int_inv_24} we obtain
\begin{equation}
{\cal I}(\Gamma_0) = i \frac{\sigma}{\rho} \frac{Q_0}{k_{\alpha\beta}^2 k_{\alpha\beta}'^2}  \hspace{40pt}  {\cal I}(\Gamma_+) = \frac{\sigma}{\rho} \frac{Q_+}{k_{\alpha\beta}^2} \hspace{40pt} {\cal I}(\Gamma_-) = \frac{\sigma}{\rho} \frac{Q_-}{k_{\alpha\beta}'^2}
\end{equation}
with
\begin{equation}
\hspace{-109pt} Q_{\pm}^2 =  (\cos \delta_{\alpha\gamma} \pm \cos \delta_{\beta\gamma})^2
\end{equation}
and 
\begin{equation}
Q_0^2\, = 
\begin{array}{|ccc|}
1 & \!\!\!\cos \delta_{\alpha\gamma}\!\!\! & \cos \delta_{\alpha\beta} \\ [4pt]
\cos \delta_{\alpha\gamma} & 1 & \cos \delta_{\beta\gamma} \\ [4pt]
\cos \delta_{\alpha\beta} & \!\!\!\cos \delta_{\beta\gamma}\!\!\! & 1 \\
\end{array}
= 36 \times ({\rm Vol}_{\, OABC})^{2} 
\end{equation}
$A$, $B$ and $C$ are the points on the round sphere corresponding to the roots $\alpha$, $\beta$ and $\gamma$, $O$ is the center of the sphere and ${\rm Vol}_{\, OABC}$ is the volume of the tetrahedron $OABC$. The vanishing of $Q_0^2$ is the necessary and sufficient condition for the three points $A$, $B$ and $C$ to lie on the same geodesic circle ({\it i.e.}, to be colinear in the sense of projective geometry).

Alternatively, define the invariants
\begin{eqnarray}
g_{\rho \sigma^2} &\hspace{-7pt} = \hspace{-7pt}& 4 \times \frac{1}{4!}\, \eta_{(AB}\, \eta_{CD)}\, \eta^{ABCD} \\
g_{\rho^2 \sigma^2} &\hspace{-7pt} = \hspace{-7pt}& 24 \times \frac{1}{4!}\, \eta_{(AB}\, \eta_{CD)}\, \eta^{ABEF} \eta_{EF}{}^{CD} 
\end{eqnarray}
The numerical factors are chosen for convenience. Explicitly, they take the nondescript forms
\begin{eqnarray}
g_{\rho\sigma^2} &\hspace{-7pt} = \hspace{-7pt}& \frac{2}{3}x_2(x_1^2-2|z_1|^2) + 4z_2\bar{z}_1^2 + 4\bar{z}_2z_1^2 + 2v_2\bar{z}_1x_1 + 2 \bar{v}_2z_1x_1 \\
g_{\rho^2\sigma^2} &\hspace{-7pt} = \hspace{-7pt}& (8|z_2|^2-|v_2|^2-\frac{2}{3}x_2^2)(x_1^2-2|z_1|^2) -12z_2\bar{v}_2\bar{z}_1x_1 - 12 \bar{z}_2v_2z_1x_1\nonumber \\ [3pt]
&\hspace{-7pt}  + \hspace{-7pt}& 8z_2x_2\bar{z}_1^2 + 8\bar{z}_2x_2z_1^2 - 2v_2x_2\bar{z}_1x_1 - 2\bar{v}_2x_2z_1x_1 - 3v_2^2\bar{z}_1^2 - 3\bar{v}_2^2z_1^2
\end{eqnarray}
The diagram corresponding to the first spherical invariant is \\ [-7pt]
\begin{figure}[H]
\centering
\scalebox{0.4}{\includegraphics{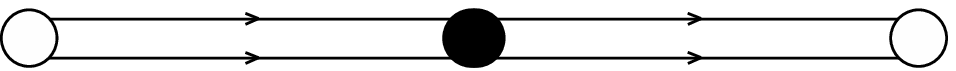}}
\put(-126,1){$4\times$}
\end{figure} 
\vspace{2pt}
\noindent whereas the second one can be represented for instance by \\ [-4pt]
\begin{figure}[H]
\centering
\scalebox{0.4}{\includegraphics{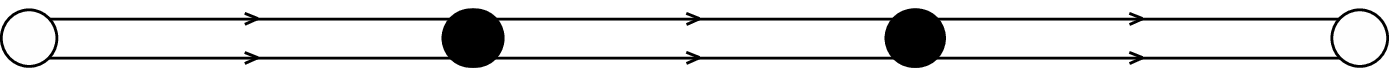}}
\put(-182,1){$24\times$}
\put(4,1){$ + \ 2 g_{\rho^2} g_{\sigma^2}$}
\end{figure} 
\noindent Any other diagram constructed from either one or both of these multiplets can be reduced to homogeneous rational polynomial expressions in terms of these basic invariants. For example,
\begin{figure}[H]
\centering
\hspace{-140pt}
\scalebox{0.4}{\includegraphics{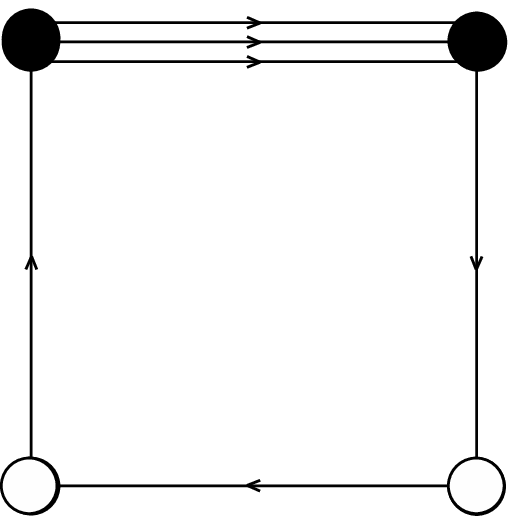}}
\put(3,27){$= -\displaystyle{\frac{1}{8} g_{\rho^2} g_{\sigma^2}}$}
\end{figure} 
\begin{figure}[H]
\centering
\hspace{-140pt}
\scalebox{0.4}{\includegraphics{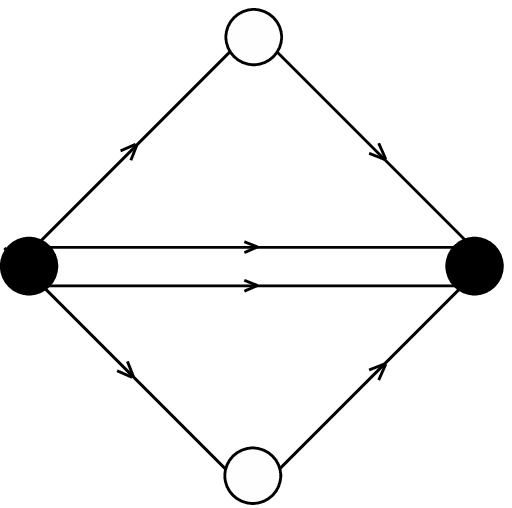}}
\put(4,25){$=\displaystyle{\frac{1}{24}(g_{\rho^2\sigma^2} + g_{\rho^2} g_{\sigma^2})}$}
\end{figure} 
\begin{figure}[H]
\centering
\hspace{-140pt}
\scalebox{0.4}{\includegraphics{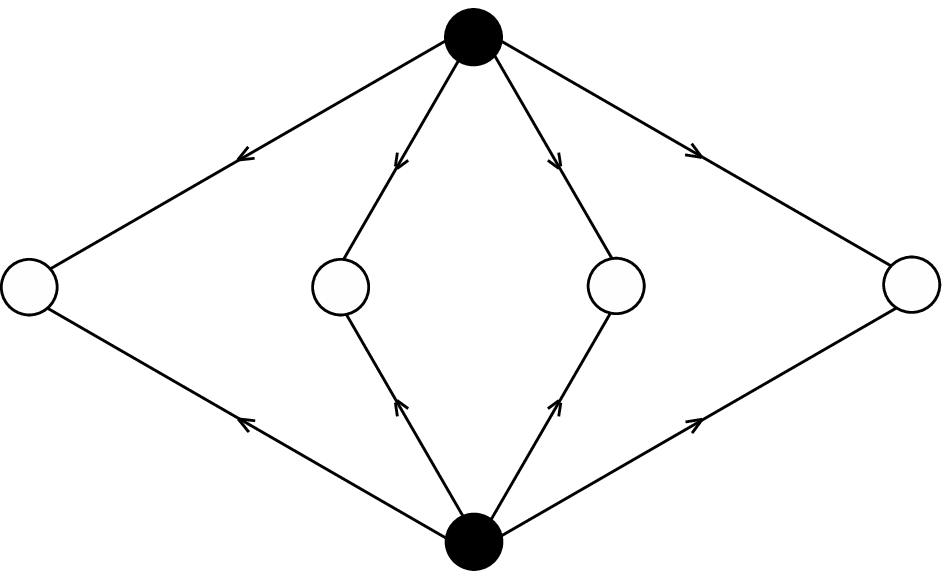}}
\put(3,30.5){$= \displaystyle{\frac{1}{96}(6g_{\rho\sigma^2}^{\,2} + 2g_{\rho^2\sigma^2} g_{\sigma^2} - g_{\rho^2} g_{\sigma^2}^{\,2})}$}
\end{figure} 
\noindent a.s.o. Note that all combinations with an odd number of ${\cal O}(2)$ vertices vanish.

To relate the two types of mixed invariants that we have introduced so far, we start from the observation that $\eta^{(4)} - \lambda\, (\eta^{(2)})^2$ is a real (w.r.t. antipodal conjugation) ${\cal O}(4)$ multiplet for any real invariant coupling scale $\lambda$. In particular, one can construct the associated ${\cal O}(4)$ basic spherical invariants 
\begin{eqnarray}
g_2(\lambda) &\hspace{-7pt} = \hspace{-7pt}& g_{\rho^2} - 3g_{\rho\sigma^2} \left(\frac{\lambda}{3}\right)  + 3g_{\sigma^2}^{\,2}\left(\frac{\lambda}{3}\right)^{\!2} \\
g_3(\lambda) &\hspace{-7pt} = \hspace{-7pt}& g_{\rho^3} - g_{\rho^2\sigma^2} \left(\frac{\lambda}{3}\right)   - 3g_{\rho\sigma^2} g_{\sigma^2}\left(\frac{\lambda}{3}\right)^{\!2} + 2g_{\sigma^2}^{\,3}\left(\frac{\lambda}{3}\right)^{\!3}
\end{eqnarray}
The following  remarkable relations hold 
\begin{eqnarray}
g_3\!\left(\frac{3e_1}{\sigma^2}\right) &\hspace{-7pt} = \hspace{-7pt}& +\frac{3}{4}\rho^2 e_1Q_-^2 \label{G_3(crit1)} \\ [2pt]
g_3\!\left(\frac{3e_2}{\sigma^2}\right) &\hspace{-7pt} = \hspace{-7pt}&  -\frac{3}{4}\rho^2 e_2Q_0^2 \\ [2pt]
g_3\!\left(\frac{3e_3}{\sigma^2}\right) &\hspace{-7pt} = \hspace{-7pt}& +\frac{3}{4}\rho^2 e_3Q_+^2 \label{G_3(crit3)}
\end{eqnarray}
with the $e_1$, $e_2$ and $e_3$ given in (\ref{W_roots}). This can be verified by expressing everything in terms of the Majorana roots and scales. The simplification that occurs at these particular couplings is quite substantial in view of the fact that $g_{\rho\sigma^2}$ and $g_{\rho^2\sigma^2}$ alone contain 36 respectively 141 terms when expressed in terms of the roots. From any two of the equations (\ref{G_3(crit1)}) through (\ref{G_3(crit3)}) one obtains the relations
\begin{eqnarray}
g_{\rho\sigma^2} &\hspace{-7pt} = \hspace{-7pt}& \rho\sigma^2 (\cos\delta_{\alpha\gamma}\cos\delta_{\beta\gamma} - \frac{1}{3}\cos\delta_{\alpha\beta}) \label{-r1s2-} \\ [1pt]
g_{\rho^2\sigma^2} &\hspace{-7pt} = \hspace{-7pt}& g_{\rho^2}g_{\sigma^2} +\frac{1}{4}\rho^2\sigma^2(Q_0^2-Q_+^2-Q_-^2) \label{-r2s2-}
\end{eqnarray}

\subsection{Rotational invariants as quantum amplitudes}

Let $| \psi_{\eta^{(2)}} \rangle$ and $| \psi_{\eta^{(4)}} \rangle$ be the spin-$1$ respectively spin-$2$ quantum coherent states associated to $\eta^{(2)}$ and $\eta^{(4)}$ according the prescription of section \ref{coherent_states}. By taking tensor products of these elementary states one can construct composite states. For example, the tensor product $| \psi_{\eta^{(2)}}\! \otimes \psi_{\eta^{(2)}} \rangle$ has a spin-$2$ component equal to $| \psi_{(\eta^{(2)})^2} \rangle$, no spin-$1$ component, and a spin-$0$ component given by $-1/(2\sqrt{3})\sigma^2 |00\rangle$. Hilbert scalar products of the quantum states formed in this way yield rotational invariants of the type discussed above:
\begin{eqnarray}
\Vert \psi_{\eta^{(2)}} \! \Vert^2 &\hspace{-7pt} = \hspace{-7pt}& \frac{1}{2} g_{\sigma^2} \label{norm-1} \\ [3pt]
\Vert \psi_{\eta^{(4)}} \! \Vert^2 &\hspace{-7pt} = \hspace{-7pt}& \frac{1}{2} g_{\rho^2} \label{norm-2} \\ [1pt]
\langle \psi_{\eta^{(4)}} | \psi_{\eta^{(4)}}\! \otimes \psi_{\eta^{(4)}} \rangle &\hspace{-7pt} = \hspace{-7pt}& \frac{9}{4\sqrt{21}} g_{\rho^3} \\
\langle \psi_{\eta^{(2)}} | \psi_{\eta^{(2)}}\! \otimes \psi_{\eta^{(4)}} \rangle &\hspace{-7pt} = \hspace{-7pt}& - \frac{3}{4\sqrt{15}} g_{\rho\sigma^2} \\
\langle \psi_{\eta^{(4)}} | \psi_{\eta^{(2)}}\! \otimes \psi_{\eta^{(2)}} \rangle &\hspace{-7pt} = \hspace{-7pt}& \frac{1}{4} g_{\rho\sigma^2} \label{<2|1x1>} \\ 
\langle \psi_{\eta^{(2)}}\! \otimes \psi_{\eta^{(2)}} | \psi_{\eta^{(4)}}\! \otimes \psi_{\eta^{(4)}} \rangle &\hspace{-7pt} = \hspace{-7pt}& \frac{1}{4\sqrt{21}} g_{\rho^2\sigma^2} - \frac{1}{4\sqrt{15}}g_{\rho^2}g_{\sigma^2} \label{<2x2|4x4>}
\end{eqnarray}
This is because these scalar products lead to expressions of the type (\ref{spherical_invs}) when written in a spherical basis. Not all such scalar products lead to independent invariants. Some vanish, yielding orthogonality relations, {\it e.g.},
\begin{eqnarray}
\langle \psi_{\eta^{(2)}} | \psi_{\eta^{(4)}} \rangle &\hspace{-7pt} = \hspace{-7pt}& 0 \\ [8pt]
\langle \psi_{\eta^{(2)}} | \psi_{\eta^{(4)}}\! \otimes \psi_{\eta^{(4)}} \rangle &\hspace{-7pt} = \hspace{-7pt}& 0 \\ [8pt]
\langle \psi_{\eta^{(4)}} | \psi_{\eta^{(2)}}\! \otimes \psi_{\eta^{(4)}} \rangle &\hspace{-7pt} = \hspace{-7pt}& 0
\end{eqnarray}
The representation of the invariants as quantum amplitudes can be put to use to derive various inequalities.  For instance, from the positive-definiteness of the Hilbert space norm it follows that
\begin{equation}
g_{\sigma^2},\ g_{\rho^2} \geq 0 \label{radial-inv}
\end{equation}
This is consistent with the conclusion that one can derive in a more direct manner by examining the Majorana-coefficient expressions of $g_{\sigma^2}$ and $g_{\rho^2}$. Other, less obvious relations follow by way of Cauchy-Schwarz inequalities on the Hilbert space. For example, from the equation (\ref{<2|1x1>}) together with (\ref{norm-1}) and (\ref{norm-2}) one gets upper and lower bounds for $g_{\rho\sigma^2}$
\begin{equation}
-\sqrt{2} \leq \frac{g_{\rho\sigma^2}}{\sqrt{g_{\rho^2}}g_{\sigma^2}}  \leq \sqrt{2} \label{angular-inv1}
\end{equation}
whereas the equation (\ref{<2x2|4x4>}) together with (\ref{norm-1}) and (\ref{norm-2}) yield upper and lower bounds for $g_{\rho^2\sigma^2}$
\begin{equation}
-\sqrt{\frac{7}{5}} (\sqrt{15}-1) \leq \frac{g_{\rho^2\sigma^2}}{g_{\rho^2}g_{\sigma^2}}  \leq \sqrt{\frac{7}{5}} (\sqrt{15}+1) \label{angular-inv2}
\end{equation}
We will henceforth refer to positive-definite invariants as in (\ref{radial-inv}) as being of {\it radial} type and to doubly-bounded invariants as in (\ref{angular-inv1}) and (\ref{angular-inv2}) as being of {\it angular} type. To underline their angular character we will sometimes use instead of the mixed invariants $g_{\rho\sigma^2}$ and $g_{\rho^2\sigma^2}$ the equivalent pair of invariants 
\begin{equation}
A = \frac{g_{\rho\sigma^2}}{\sqrt{3g_{\rho^2}}g_{\sigma^2}} \qquad\mbox{and}\qquad
B = - \frac{g_{\rho^2\sigma^2}}{3g_{\rho^2}g_{\sigma^2}}
\end{equation}
The numerical factors have been chosen for later convenience.

\section{An ${\cal O}(2)\oplus{\cal O}(2)$-based Swann bundle} \label{O2+O2_Swann}

In this section we review the generalized Legendre transform construction of the 8-dimensional Swann bundle with two abelian tri-holomorphic isometries, generated by the $F$-function
\begin{equation}
F = \frac{1}{2\pi i}\oint \frac{d\zeta}{\zeta} \frac{(\eta^{(2)}_1)^2}{\eta^{(2)}_2}
\end{equation}
depending on two ${\cal O}(2)$ multiplets 
\begin{equation}
\eta^{(2)}_I = \frac{\bar{z}_I}{\zeta} + x_I - z_I\zeta 
\end{equation}
($I=1,2$). 
The integration contour winds around the roots of $\eta^{(2)}_2$ in such a way that the integral yields a real outcome. That the resulting hyperk\"{a}hler variety has a Swann bundle structure follows from the fact that $F$ scales with weight one under the weight-one scaling of the two ${\cal O}(2)$ multiplets \cite{MA1}. This Swann bundle was used by Calderbank and Pedersen in \cite{Calderbank:2001uz} to classify  selfdual Einstein metrics with two commuting isometries and by Anguelova, Ro\v{c}ek and Vandoren in \cite{Anguelova:2004sj} to describe the geometry of the classical moduli space of the universal hypermultiplet in string theory compactifications. We retrace here the basic steps of their construction and compute,  additionally, the hyperk\"{a}hler potential of the metric.

The residue theorem yields for the contour integral the expression
\begin{equation}
F = \frac{r_1^2r_2^2- (\vec{r}_1\!\cdot\!\vec{r}_2)^2}{2r_2|z_2|^2} + \frac{r_2(z_1\bar{z}_2-z_2\bar{z}_1)^2}{|z_2|^4}
\end{equation}
where
\begin{equation}
\vec{r}_I = (z_I+\bar{z}_I, -i(z_I-\bar{z}_I),\, x_I)
\end{equation}
with $I=1,2$ are the standard $\mathbb{R}^3$ vectors associated to $\eta^{(2)}_1$ and $\eta^{(2)}_2$. For the first derivatives of $F$ with respect to $x_1$ and $x_2$ we get
\begin{eqnarray}
x_1\frac{\partial F}{\partial x_1} &\hspace{-7pt} = \hspace{-7pt}& - \frac{2r_2|z_1|^2}{|z_2|^2} + \frac{2(z_1\bar{z}_2-z_2\bar{z}_1)^2}{r_2|z_2|^2} + \frac{r_1^2r_2^2- (\vec{r}_1\!\cdot\!\vec{r}_2)^2}{2r_2|z_2|^2} + \frac{2r_1^2}{r_2} \\
x_2\frac{\partial F}{\partial x_2} &\hspace{-7pt} = \hspace{-7pt}& \phantom{+}\frac{2r_2|z_1|^2}{|z_2|^2} - \frac{2(z_1\bar{z}_2-z_2\bar{z}_1)^2}{r_2|z_2|^2} + \frac{r_2(z_1\bar{z}_2-z_2\bar{z}_1)^2}{|z_2|^4} - \frac{2(\vec{r}_1\!\cdot\!\vec{r}_2)^2}{r_2^3}
\end{eqnarray}
The hyperk\"{a}hler potential of the Swann bundle follows from the generalized Legendre transform prescription, which in this case reads
\begin{equation}
K = F - x_1(u_1+\bar{u}_1) - x_2(u_2+\bar{u}_2)
\end{equation}
with
\begin{eqnarray}
\frac{\partial F}{\partial x_1} &\hspace{-7pt} = \hspace{-7pt}& u_1 + \bar{u}_1 \\
\frac{\partial F}{\partial x_2} &\hspace{-7pt} = \hspace{-7pt}& u_2 + \bar{u}_2
\end{eqnarray}
Assemblying all the pieces above, we obtain
\begin{equation}
K = - \frac{2(\vec{r}_1\times\vec{r}_2)^2}{r_2^3} \label{HK-pot_O2+O2}
\end{equation}
The dependence on the holomorphic coordinates is implicit, the $SO(3)$ invariance, on the other hand, is manifest.

The metric components in the holomorphic coordinate basis $\{z_I,u_I\}_{I=1,2}$ are related to the second derivatives of $F$. This coordinate basis obscures the $SO(3)$ structure as well as the abelian isometries induced by the ${\cal O}(2)$ multiplets. By switching instead to the real coordinate basis provided by $\vec{r}_I$ and $\psi_I = \mbox{Im}\, u_I$, the symmetries become more transparent, at the expense of the holomorphic structure. In this new basis, the metric takes the following generalized Gibbons-Hawking form
\begin{equation}
ds^2 = \Phi_{IJ}\, d\vec{r}_I \!\cdot\! d\vec{r}_J + (\Phi^{-1})_{IJ}(d\psi_I+\vec{A}_{IK} \!\cdot\! d\vec{r}_K)(d\psi_J+\vec{A}_{JL} \!\cdot\! d\vec{r}_L) \label{Gibbons-Hawking}
\end{equation}
with the generalized Bogomol'nyi conditions
\begin{equation}
\vec{\nabla}_I \times \vec{A}_{KJ} = - \vec{\nabla}_I \Phi_{KJ} \qquad\mbox{and}\qquad \vec{\nabla}_I \Phi_{KJ} = \vec{\nabla}_J \Phi_{KI}
\end{equation}
The operators $\vec{\nabla}_I = \partial/\partial \vec{r}_I$ are usual $\mathbb{R}^3$ gradients. This form holds generically for any generalized Legendre transform construction based exclusively on ${\cal O}(2)$ multiplets. In our particular case, the Higgs-field matrix takes the form
\begin{equation}
(\Phi_{IJ}) = \frac{2}{r_2} \!
\left( \!
\begin{array}{cc}
-1 & \displaystyle{\frac{\vec{r}_1\!\cdot\!\vec{r}_2}{r_2^2}} \\ [10pt]
\displaystyle{\frac{\vec{r}_1\!\cdot\!\vec{r}_2}{r_2^2}} & \displaystyle{\frac{r_1^2r_2^2-3 (\vec{r}_1\!\cdot\!\vec{r}_2)^2}{2r_2^4}}
\end{array}
\!\! \right)
\end{equation}

\section{An ${\cal O}(2)\oplus{\cal O}(4)$-based Swann bundle} \label{O4+O4_Swann}

\subsection{The hyperk\"{a}hler potential}

In \cite{Anguelova:2004sj} it was conjectured, based on a symmetry argument, that the nonperturbative universal hypermultiplet moduli space metric due to five-brane instantons is a certain deformation of the hyperk\"{a}ler variety that is generated, through the generalized Legendre transform, by the $F$-function
\begin{equation}
F =  \oint_{\Gamma} \frac{d\zeta}{\zeta} \frac{(\eta^{(2)})^2}{\sqrt{\eta^{(4)}}} \label{F_UH}
\end{equation}
The contour $\Gamma$ around the branch-cuts of $\sqrt{\eta^{(4)}}$ is chosen in such a way that the outcome of the contour integration is real. Since $F$ scales with weight $1$ under the scaling transformation 
\begin{equation}
\eta^{(2)} \longrightarrow \lambda \eta^{(2)} \qquad\qquad \eta^{(4)} \longrightarrow \lambda^2 \eta^{(4)} \label{scale24}
\end{equation}
the resulting $8$-dimensional hyperk\"{a}hler variety will have a Swann bundle structure \cite{MA1}. Swann bundles possess a so-called hyperk\"{a}hler  potential, a function defined up to the addition of a constant which is simultaneously a K\"{a}hler potential for each complex structure compatible with the hyperk\"{a}hler  structure.  For Swann bundles, the generalized Legendre transform construction produces the hyperk\"{a}hler potential.  The generalized Legendre transform relations for an ${\cal O}(2)\oplus{\cal O}(4)$ theory read
\begin{equation}
K = F - u_2v_2 - \bar{u}_2\bar{v}_2 - (u_1+\bar{u}_1)x_1
\end{equation}
with
\begin{eqnarray}
\frac{\partial F}{\partial x_1} &\hspace{-7pt} = \hspace{-7pt}& u_1+\bar{u}_1 \label{Leg-rel_UH1} \\
\frac{\partial F}{\partial v_2} &\hspace{-7pt} = \hspace{-7pt}& u_2 \label{Leg-rel_UH2} \\
\frac{\partial F}{\partial x_2} &\hspace{-7pt} = \hspace{-7pt}& 0 \label{dFdx=0}
\end{eqnarray}
The holomorphic coordinates are $z_1$, $u_1$, $z_2$, $u_2$. We differentiate by means of an index 1 or 2 between quantities related to the ${\cal O}(2)$ and the ${\cal O}(4)$ multiplet respectively, and use in general the notations established in sections \ref{O2-invs}, \ref{O4-invs}  and \ref{mixed-invs}.

To evaluate $F$, observe that we can write
\begin{equation}
2F = z_1^2F_{z_1z_1} + 2z_1x_1F_{z_1x_1} + (x_1^2-2|z_1|^2)F_{x_1x_1} + 2\bar{z}_1x_1F_{\bar{z}_1x_1}+\bar{z}_1^2F_{\bar{z}_1\bar{z}_1} 
\end{equation}
The double derivatives of $F$ can in turn be further expressed as follows
\begin{equation}
F_{z_1z_1} = 4{\cal I}^{(1)}_2 \qquad\qquad 
F_{z_1x_1} = - 4{\cal I}^{(1)}_1 \qquad\qquad
F_{x_1x_1} = 4{\cal I}^{(1)}_0 
\end{equation}
in terms of the purely ${\cal O}(4)$ elliptic integrals
\begin{equation}
{\cal I}^{(1)}_m = \oint_{\Gamma} \frac{d\zeta}{\zeta} \frac{\zeta^m}{2\sqrt{\eta^{(4)}(\zeta)}}
\end{equation}
These integrals have been discussed in detail in \cite{MM3}, where we evaluated them in terms of Weierstrass elliptic integrals. We quote here the results
\begin{eqnarray}
{\cal I}_0^{(1)} &\hspace{-7pt} = \hspace{-7pt}& 2\omega \label{I0-c} \\ [8pt]
{\cal I}_1^{(1)} &\hspace{-7pt} = \hspace{-7pt}& \frac{\pi(x_+)+\pi(x_-)}{\sqrt{z}}  \label{I1-c} \\
{\cal I}_2^{(1)} &\hspace{-7pt} = \hspace{-7pt}& -\frac{2\eta + (x_+\!+x_-)\omega - (v_-\!+iv_+)[\pi(x_+)+\pi(x_-)]}{2z}  \label{I2-c}
\end{eqnarray}
$x_{\pm}$ and $v_{\pm}$ are related to the Majorana coefficients of $\eta^{(4)}$, 
\begin{equation}
x_{\pm} = \frac{x_2\pm6|z_2|}{3} \qquad\qquad v_+ = \mbox{Im}\frac{v_2}{\sqrt{z_2}} \qquad\qquad v_- = \mbox{Re}\frac{v_2}{\sqrt{z_2}} \label{xpm_vpm}
\end{equation}
The $\pi$-function is the Weierstrass representation of Jacobi's version of the elliptic integral of third kind. We defined it and studied its properties in \cite{MM1}, to which we refer for further details. Putting things together, we eventually obtain that
\begin{eqnarray}
F &\hspace{-7pt} = \hspace{-7pt}& 4(z_{1+}^2\! - z_{1-}^2)\, \eta + 4(x_1^2+x_-z_{1+}^2\! - x_+z_{1-}^2)\, \omega \nonumber \\ [2pt]
 &\hspace{-7pt} + \hspace{-7pt}& 2\, \mbox{Re}\! \left[ \frac{z_1}{\sqrt{z_2}} \left(\frac{v_2z_1}{z_2}-4x_1\right) \!\right]\! \pi(x_-) + 2i\, \mbox{Im}\! \left[ \frac{z_1}{\sqrt{z}_2} \left(\frac{v_2z_1}{z_2}-4x_1\right) \!\right]\!  \pi(x_+) \label{F'_UH}
\end{eqnarray}
where we introduced the additional notations
\begin{equation}
z_{1+} = \mbox{Im}\frac{z_1}{\sqrt{z_2}} \qquad\qquad z_{1-} = \mbox{Re}\frac{z_1}{\sqrt{z_2}} \label{z1pm}
\end{equation}

Then, similarly to the Atiyah-Hitchin case, by means of the elliptic differentiation formulas given in the Appendix, we compute the following derivatives of (\ref{F'_UH})
\begin{eqnarray}
\frac{\partial F}{\partial x_1} &\hspace{-7pt} = \hspace{-7pt}& 8 [x_1\omega - z_{1-}\pi(x_-) - iz_{1+}\pi(x_+)] \\ [3pt]
\frac{\partial F}{\partial v_2} &\hspace{-7pt} = \hspace{-7pt}& \frac{ 2M \eta^* + 2N \omega^* 
+ (z_{1-}\!+iz_{1+})^2[\pi(x_-)+\pi(x_+)] }{\sqrt{z_2}}  \\
\frac{\partial F}{\partial x_2} &\hspace{-7pt} = \hspace{-7pt}&  - 2 (g_{\rho\sigma^2}\, \eta^* - g_{\rho^2\sigma^2}\, \omega^*)  \label{F'x_UH}
\end{eqnarray}
where $\eta^*$ and $\omega^*$ are defined in~(\ref{omet1}) and $M=M_-\!+iM_+$, $N=N_-\!+iN_+$, with
\begin{eqnarray*}
M_{\pm} &\hspace{-7pt} = \hspace{-7pt}& \phantom{+}\frac{(2g_{\rho^2}x_{\pm}+3g_{\rho^3})[g_{\sigma^2} - 3(x_+\!-x_-)z_{1\pm}^2] - 3x_{\pm}^2g_{\rho\sigma^2} - x_{\pm}g_{\rho^2\sigma^2}}{3(x_+\!-x_-)v_{\pm}}  \\ [4pt]
N_{\pm} &\hspace{-7pt} = \hspace{-7pt}& \frac{-(9g_{\rho^3}x_{\pm}+2g_{\rho^2}^2)[g_{\sigma^2}-3(x_+\!-x_-)z_{1\pm}^2] + (3g_{\rho^2}x_{\pm}+9g_{\rho^3})g_{\rho\sigma^2} + (3x_{\pm}^2-2g_{\rho^2})g_{\rho^2\sigma^2}}{3(x_+\!-x_-)v_{\pm}}
\end{eqnarray*}
The imaginary parts $i M_+$ and $i N_+$ of the coefficients $M$ and $N$ are conjugates of the corresponding real parts $M_-$ and $N_-$ under the $\mathbb{Z}_2$ action 
\begin{eqnarray}
x_-	& \longleftrightarrow & x_+ \nonumber \\
v_-	& \longleftrightarrow & iv_+ \nonumber \\
z_{1-} &\longleftrightarrow & iz_{1+} 
\end{eqnarray}
The r.h.s. of equation (\ref{F'x_UH}) is manifestly $SO(3)$-invariant. That this should be so can be argued independently, without resorting to direct calculation, as follows: commuting the derivative with the integral, one obtains the integral representation
\begin{equation}
\frac{\partial F}{\partial x_2} = -\frac{1}{2} \oint_{\Gamma} \frac{d\zeta}{\zeta} \frac{(\eta^{(2)})^2}{(\eta^{(4)})^{3/2}}
\end{equation}
Under the scaling transformation (\ref{scale24}) this integral transforms with weight $-1$. According to the discussion following equation (\ref{Penrs-int}), it should then result in a $SO(3)$-invariant quantity.  The equation (\ref{F'x_UH}) also provides us with a good opportunity to advertize the superiority of the Weierstrass approach. Had we expressed the multiplets in terms of the Majorana roots and evaluated the derivative of $F$ with respect to $x_2$ within the Legendre frame we would have obtained an expression with 709 terms!

Note the structural similarity between the Jacobi terms in equation (\ref{F'_UH}) and the corresponding equation in the Atiyah-Hitchin case \cite{MM1}. The same mechanism gives us now the hyperk\"{a}hler potential: the Jacobi terms cancel against the quadratic terms in the Legendre transform when the Legendre relations (\ref{Leg-rel_UH1}) and (\ref{Leg-rel_UH2}) are used.  The resulting K\"{a}hler potential is
\begin{equation}
K = -\frac{4}{3}(g_{\rho^2\sigma^2}+4g_{\rho^2}g_{\sigma^2})\eta^* + 4(g_{\rho^2}g_{\rho\sigma^2}+6g_{\rho^3}g_{\sigma^2}) \omega^* - 4(x_+\!+x_-)( g_{\rho\sigma^2}\, \eta^* - g_{\rho^2\sigma^2}\omega^*) \label{K-mixed}
\end{equation}
On the other hand, the generalized Legendre transform relation (\ref{dFdx=0}) reads
\begin{equation}
g_{\rho\sigma^2}\eta^* = g_{\rho^2\sigma^2}\omega^*   \label{aux_UH}
\end{equation}
Upon resorting to it, the $(x_+\!+x_-)$-dependent terms in (\ref{K-mixed}) drop out and the resulting hyperk\"{a}hler potential takes the remarkably compact manifestly $SO(3)$-invariant form
\begin{equation}
K = -\frac{4}{3}(g_{\rho^2\sigma^2}+4g_{\rho^2}g_{\sigma^2})\eta^* + 4(g_{\rho^2}g_{\rho\sigma^2}+6g_{\rho^3}g_{\sigma^2}) \omega^* \label{HK-pot-UH}
\end{equation}
That the hyperk\"{a}hler potential must be invariant under $SO(3)$ transformations can be argued on general grounds, and this provides an additional validation for our result.

The hyperk\"{a}hler holomorphic $(2,0)$-form takes the Darboux form
\begin{equation}
\omega^+ = dz_1\wedge du_1 + dz_2\wedge du_2
\end{equation}
Just as in the Atiyah-Hitchin case it is worthwhile to perform the following holomorphic symplectomorphism
\begin{equation}
U_2 = u_2\sqrt{z_2} \qquad\quad Z_2 = 2\sqrt{z_2} 
\end{equation}
In the new holomorphic coordinate basis,
\begin{equation}
\omega^+ = dz_1\wedge du_1 + dZ_2\wedge dU_2
\end{equation}
and the conformal homothetic Killing vector field reads
\begin{equation}
X = 2\left( z_1\frac{\partial}{\partial z_1} + \bar{z}_1\frac{\partial}{\partial \bar{z}_1} + Z_2\frac{\partial}{\partial Z_2} + \bar{Z}_2\frac{\partial}{\partial \bar{Z}_2} \right)
\end{equation}
One can check explicitly that the K\"{a}hler potential $K$ is an eigenfunction of $X$, {\it i.e.},
\begin{equation}
X(K) = 2K
\end{equation}
Besides the Swann bundle structure, the variety has an additional abelian tri-holomorphic isometry that is due to the presence of the ${\cal O}(2)$ multiplet. This is generated by the Killing vector field
\begin{equation}
\tilde{X} = i \left( \frac{\partial}{\partial u_1} - \frac{\partial}{\partial \bar{u}_1} \right)
\end{equation}

One can go further and compute the metric explicitly in this holomorphic coordinate basis. For that, one needs to compute the second derivatives of $F$ with respect to the Majorana coefficients of the two multiplets. In principle, this is a straightforward task, since all necessary tools have already been developed \cite{MM1}. Unfortunately we have not been able to cast the result in a presentable compact form. A reasonable guess is that, nevertheless, such a form is very likely to exist, perhaps in  a coordinate basis better adapted to the many symmetries of the problem than our own.

\subsection{Asymptotic expansions}

The single ${\cal O}(2)$ invariant $\sigma=r_1$ is of radial type. On the ${\cal O}(4)$ side, there are two  invariants of radial type, namely $r_2$ and $r'_2$ defined in (\ref{r&r'}). Other ${\cal O}(4)$ invariants such as $g_{\rho^2}$, $g_{\rho^3}$ and $\eta$ have theta-function representations which allow one to express them in terms of  $r_2$ and $r'_2$ in the form of infinite Lambert-type series. By constrast, the mixed invariants are essentially of angular type. The following table summarizes the various radial and angular invariants associated to an ${\cal O}(2)\oplus{\cal O}(4)$ system of multiplets 
\vspace{4pt}
\begin{center}
\begin{tabular}{|c|c|c|} 
${\cal O}(2)$ \mbox{invariant\phantom{s}} & ${\cal O}(4)$ \mbox{invariants} & \mbox{mixed invariants} \\ \hline\hline
$r_1$ & $r_2$, $r'_2$ & $A$, $B$ \\ \hline
\multicolumn{2}{|c|}{\mbox{radial}} & \mbox{angular} 
\end{tabular}
\end{center} 
\vspace*{10pt}
We want to investigate the behavior of the generalized Legendre transform equation (\ref{aux_UH}) and the hyperk\"{a}hler potential (\ref{HK-pot-UH}) in and around the asymptotic limits $r_2>>r'_2$ and $r_2<<r'_2$. This is facilitated in a decisive manner by their manifest $SO(3)$ invariance. In practice, the two asymptotic regions are probed by expanding the ${\cal O}(4)$ invariants in the {\it nome} $q$ respectively the complementary {\it nome} $q'$, as explained at the end of section \ref{O4-invs}. Specifically, the $q$-series expansions for $g_{\rho^2}$ and $g_{\rho^3}$ are given by (\ref{g2-Lambert}) and (\ref{g3-Lambert}), while for $\eta$ we have (see {\it e.g.} \cite{MR1054205})
\begin{equation}
\eta = \frac{\pi^2}{12\omega} \left(1 - 24 \sum_{n=1}^{\infty} n \frac{q^{2n}}{1-q^{2n}} \right) \label{eta-Lambert}
\end{equation}
To obtain the $q'$-series expansions for the Weierstrass coefficients one can use the fact that they are invariant unde the modular transformation $\tau' = -1/\tau$ and so the equations (\ref{g2-Lambert}) and (\ref{g3-Lambert}) still hold if one replaces $\omega$ and $q$ with $\omega'$ and $q'$. Similarly, the equation (\ref{eta-Lambert}) still holds if one replaces $\eta$, $\omega$ and $q$ by $\eta'$, $\omega'$ and $q'$. This yields a $q'$-series expansion for $\eta'$. Furthermore, $\eta'$ is related to $\eta$ by means of the Legendre identity
\begin{equation}
\begin{array}{|cc|}
\omega'	& \omega \\
\eta'		& \eta 
\end{array} 
=i\frac{\pi}{2} \label{Legendre-id}
\end{equation}
 which then allows us to write down a $q'$-series expansion for the latter. Clearly, one can perform these expansions virtually to any order.

The $q'$-series expansion of equation (\ref{aux_UH}) yields
\begin{equation}
\frac{B}{A} = 1 - \frac{288(3r_2-r'_2)}{r_2}e^{-2r'_2/r_2} + \frac{6912(39r_2^2-26r_2r'_2+5r_2^{\prime 2})}{r_2^2}e^{-4r'_2/r_2} + \cdots \label{B-q'-exp}
\end{equation}
The limit $q' \rightarrow 0$ corresponds to the pinching of the $b$-cycle of the torus associated to the ${\cal O}(4)$ multiplet. In terms of the roots of $\eta^{(4)}$ this limit corresponds to $\alpha \rightarrow \beta$, while the ${\cal O}(4)$ multiplet degenerates into the square of an ${\cal O}(2)$ multiplet. Putting $\alpha = \beta$ in equation (\ref{-r1s2-}) and letting $\delta = \delta_{\alpha\gamma} = \delta_{\beta\gamma}$ be the Fubini-Study distance on the Riemann sphere between the confounding limit point and the $\eta^{(2)}$ root $\gamma$, we get
\begin{equation}
A = \cos^2\!\delta - \frac{1}{3} \label{A-zero-ord}
\end{equation}
Doing the same in equation (\ref{-r2s2-}) we obtain that $B=A$, in agreement with the zero-order term in the expansion (\ref{B-q'-exp}).

On the other hand, solving equation (\ref{aux_UH}) for $g_{\rho^2\sigma^2}$, substituting the result in the formula (\ref{HK-pot-UH}) for the hyperk\"{a}hler potential and then performing a $q'$-series expansion, we get
\begin{eqnarray}
K &\hspace{-7pt} = \hspace{-7pt}& 2\!\left(\! A-\frac{2}{3}\right)\!\frac{r_1^2}{r_2} - 6A\frac{r_1^2}{r'_2} \nonumber \\ [2pt]
&\hspace{-7pt} + \hspace{-7pt}& A\frac{r_1^2}{r'_2} \frac{144(5r_2^2 - 7 r_2r'_2 + 2r_2^{\prime 2})}{r_2^2} e^{-2r'_2/r_2} \nonumber \\
&\hspace{-7pt} - \hspace{-7pt}& A\frac{r_1^2}{r'_2}\frac{432(285r_2^3 - 678r_2^2r'_2 + 416r_2r_2^{\prime 2} - 80r_2^{\prime 3})}{r_2^3}e^{-4r'_2/r_2} + \cdots
\end{eqnarray}
When $r_2<<r'_2$, the dominating contribution comes from the non-exponential term. This, in turn,  contains a leading and a sub-leading part. Using the zero-order result (\ref{A-zero-ord}), the leading part of the hyperk\"{a}hler potential can be cast in the form 
\begin{equation}
K_0 = 2\!\left(\! A-\frac{2}{3}\right)\!\frac{r_1^2}{r_2} = -\frac{2r_1^2 \sin^2\!\delta}{r_2}
\end{equation}
Observe that this coincides precisely with the hyperk\"{a}hler potential (\ref{HK-pot_O2+O2}) of the ${\cal O}(2)\oplus{\cal O}(2)$ model discussed in section \ref{O2+O2_Swann}!

Let us now look at the other asymptotic region. The $q$-series expansion of equation (\ref{aux_UH}) yields
\begin{equation}
\frac{B}{A} = \frac{7}{5} - \frac{504}{5}e^{-2\pi^2r_2/r'_2} + \frac{101808}{5} e^{-4\pi^2r_2/r'_2} + \cdots \label{B-q-exp}
\end{equation}
while the expansion of the hyperk\"{a}hler potential gives
\begin{equation}
K = \frac{2}{5}\!\left(\! A-\frac{10}{3}\right) \!\frac{r_1^2}{r_2} + A \frac{r_1^2}{r_2} \frac{216}{5} e^{-2\pi^2r_2/r'_2} -  A  \frac{r_1^2}{r_2} \frac{14832}{5} e^{-4\pi^2r_2/r'_2} + \cdots
\end{equation}
In terms of the roots of $\eta^{(4)}$ the limit $q \rightarrow 0$ corresponds to $\alpha \rightarrow -1/\bar{\beta}$, while the ${\cal O}(4)$ multiplet degenerates into {\it minus} the square of an ${\cal O}(2)$ multiplet. Putting $\alpha=-1/\bar{\beta}$ and using that $\delta_{-1/\bar{\beta}\,\gamma} = \pi - \delta_{\beta\gamma}$ in the equations (\ref{-r1s2-}) and (\ref{-r2s2-}), we get that $B=-A$, which seems to be in contradiction to the leading term of (\ref{B-q-exp}). The resolution of this paradox comes from realizing that while on one hand $B=-A$ is a purely zero-order result, no corrections whatsoever being taken into account during its derivation, on the other hand the leading term in (\ref{B-q-exp}) is fundamentally a first-order result in $q^2$. Indeed, we have the $q$-series expansions
\begin{eqnarray}
\eta^*\Delta &\hspace{-7pt} = \hspace{-7pt}& \hspace{3pt} 112\pi\left(\frac{\pi}{2\omega}\right)^7 [q^2+66q^4+\cdots] \\ [4pt]
\omega^*\Delta &\hspace{-7pt} = \hspace{-7pt}& -80\pi\left(\frac{\pi}{2\omega}\right)^5 [q^2+18q^4+\cdots]
\end{eqnarray}
There are no zero-order terms to begin with. The leading term in (\ref{B-q-exp}) follows from substituting these expansions in equation (\ref{aux_UH}) and truncating consistently to first-order in $q^2$. This is to be contrasted with the situation at the other asymptotic region, where we have the $q'$-series expansions
\begin{eqnarray}
\eta^*\Delta &\hspace{-7pt} = \hspace{-7pt}&  -\frac{2i}{3} \left(\frac{\pi}{2\omega'}\right)^7[1-168(\ln q'+3)q^{\prime 2}+\cdots] \\ [4pt]
\omega^*\Delta &\hspace{-7pt} = \hspace{-7pt}& -\frac{2i}{3} \left(\frac{\pi}{2\omega'}\right)^5 [1+120(\ln q'+2)q^{\prime 2}+\cdots]
\end{eqnarray}
which do have zero-order terms in $q'^2$  and where the resulting leading term of (\ref{B-q'-exp}) is of truly zero-order nature.

\section{APPENDIX: Elliptic differentation formulas}

In these notes, we use the following form for the Weierstrass cubic
\begin{equation}
Y^2 = X^3 - g_2X - g_3
\end{equation}
The role of the complete elliptic integrals in the Weierstrass formalism is played by
\begin{eqnarray}
\omega &\hspace{-7pt} = \hspace{-7pt}& \phantom{+}\int_{e_2}^{e_3} \frac{dX}{2Y} \\
\eta &\hspace{-7pt} = \hspace{-7pt}& -\int_{e_2}^{e_3} X\frac{dX}{2Y} \\
\pi(X_0) &\hspace{-7pt} = \hspace{-7pt}& -\int_{e_2}^{e_3} \frac{Y_0}{X-X_0}\frac{dX}{2Y} = 
\begin{array}{|cc|}
u_0		& \omega \\
\zeta(u_0)	\!&\! \zeta(\omega) 
\end{array} 
\end{eqnarray}
where $u_0$ is the image of $(X_0,Y_0)$ through the Abel-Jacobi map. The notation $\pi(X_0)$ conceals a sign ambiguity. To avoid that, one should write $\pi(X_0,Y_0)$ or $\pi(u_0)$ instead. The periods $\omega$ and $\eta$ are functions of the Weierstrass coefficients $g_2$ and $g_3$ while $\pi(X)$ is additionally a function of $X$. In \cite{MM1} we proved the following differentiation formulas 
\begin{eqnarray}
d\omega 	&\hspace{-7pt} = \hspace{-7pt}& \frac{1}{2} [3\,\omega^*dg_3 - \eta^*dg_2] \\ [3pt]
d\eta 	&\hspace{-7pt} = \hspace{-7pt}& \frac{1}{2} [\eta^*dg_3 - g_2\omega^*dg_2]
\end{eqnarray}
and
\begin{equation}
d\pi(X) = \frac{(\eta+X\omega)}{2Y}dX + \frac{(g_2X+3g_3)\omega^*-X^2\eta^*}{2Y}dg_2 + \frac{(3X^2-2g_2)\omega^*-X\eta^*}{2Y}dg_3
\end{equation}
We use the shorthand notations
\begin{equation}
\eta^* = \frac{2g_2^2\omega-9 g_3\eta}{\Delta} \qquad\mbox{and}\qquad \omega^* = \frac{3g_3\omega-2g_2\eta}{\Delta} \label{omet1}
\end{equation} 
where $\Delta = 4g_2^3 - 27 g_3^2$ is the Weierstrass discriminant.

\vskip20pt
\noindent {\large \bf Acknowledgements} \\ [10pt]
It is my pleasure to thank Martin Ro\v{c}ek for valuable suggestions, enlightening discussions, encouragement and support.

\bibliographystyle{utphys}
\bibliography{article2}

\end{document}